\documentclass[11pt]{amsart}
\usepackage{amssymb}
\usepackage{textcomp}

\newtheorem{theorem}{Theorem}
\newtheorem*{theorem*}{Theorem}
  \swapnumbers
  \theoremstyle{remark}
\newtheorem{remark}[subsection]{Remark}
\newtheorem*{definition}{Definition}
\newtheorem{taller}[subsection]{}
\newenvironment{lemma}{\begin{taller}{\bf Lemma.}\em}{\end{taller}}

\newenvironment{blanko}[1]{\begin{taller}{\bf #1}\rm}{\end{taller}}

\newenvironment{blanko*}[1]{%
\begin{list}{\bf {#1} }{%
\setlength{\labelsep}{0mm}\setlength{\leftmargin}{0mm}%
\setlength{\labelwidth}{0mm}\setlength{\listparindent}{\parindent}%
\setlength{\parsep}{\parskip}\setlength{\partopsep}{0mm}}%
\item}{\end{list}}

\newcommand{\isopil}{\stackrel{\raisebox{0.1ex}[0ex][0ex]{\(\sim\)}}%
                        {\raisebox{-0.15ex}[0.28ex]{\(\rightarrow\)}}}

\newcommand{\lowerstar}{_{\raisebox{-0.33ex}[-0.5ex][0ex]{\(\ast\)}}}
\newcommand{\df}{\: {\raisebox{0.255ex}{\normalfont\scriptsize :\!\!}}=}

\newcommand{\un}{\underline}

\providecommand{\kat}[1]{\text{\textbf{\textsl{#1}}}}

\newcommand{\into}{\hookrightarrow}

\newcommand{\op}{^{\text{{\rm{op}}}}}
\newcommand{\Hom}{\operatorname{Hom}}
\newcommand{\Hocolim}{\operatorname{hocolim}}
\newcommand{\id}{\operatorname{id}}
\newcommand{\End}{\operatorname{End}}

\newcommand{\Mon}{\kat{Mon}}
\newcommand{\Top}{\kat{Top}}
\newcommand{\sSet}{\kat{sSet}}
\newcommand{\Cat}{\kat{Cat}}
\newcommand{\sCat}{\kat{sCat}}
\newcommand{\CCat}{\kat{CCat}}
\newcommand{\Ho}{\text{{\rm{Ho}}}}

\newcommand{\ABim}{\kat{Bimod}_A}

\newcommand{\REnd}{\R\un\End}
\newcommand{\RHom}{\R\un\Hom}

\newcommand{\M}{M} 
\newcommand{\tensor}{\otimes} 
\newcommand{\unit}{I} 

\renewcommand{\Mc}{\M^{\text{c}}} 
\newcommand{\Mci}{M^{\text{\rm \textcent}}} 

\newcommand{\R}{\mathbb{R}}
\newcommand{\HH}{H\! H}

\begin{document}

\title[Simplicial localization of monoidal structures]%
{Simplicial localization of monoidal structures, \\
and a non-linear version of Deligne's conjecture}

\author{Joachim Kock}
\author{Bertrand To\"en}

\address{
Laboratoire
J.~A.~Dieudonn\'e\\
Universit\'e de Nice Sophia-Antipolis\\
Parc Valrose, 06108 Nice c\'edex 2, France.}

\email{kock@math.unice.fr} \email{toen@math.unice.fr}

%

\begin{abstract}
  We show that if $(\M,\tensor,\unit)$ is a monoidal model category
  then $\REnd_{\M}(\unit)$ is a (weak) $2$-monoid in $\sSet$.  This
  applies in particular when $\M$ is the category of $A$-bimodules
  over a simplicial monoid $A$: the derived endomorphisms of $A$ then
  form its Hochschild cohomology, which therefore becomes a simplicial
  \mbox{$2$-monoid}.
\end{abstract}

\maketitle

\vspace*{-1ex}

\setcounter{section}{-1}
\section{Introduction}

\begin{blanko*}{Deligne's conjecture.}
Deligne's conjecture (stated informally in a letter in 1993)
states that the Hochschild cohomology $\HH(A)$ of an associative
algebra $A$ is a $2$-algebra --- this means that up to homotopy it
has two compatible multiplication laws.

Various versions of this conjecture have been proved,
cf.~e.g.~\cite{Getzler-Jones:9403}, \cite{Gerstenhaber-Voronov:9409},
\cite{Tamarkin:9803}, \cite{Kontsevich:9904}, \cite{Voronov:9908},
\cite{McClure-Smith:9910}, \cite{Kontsevich-Soibelman:0001},
\cite{Hu-Kriz-Voronov}, \cite{Berger-Fresse:0109}.  All these proofs
are technical, and a more conceptual proof would certainly be
desirable (we refer for example to \cite{Batanin:0207} for a
conceptual point of view on Deligne's conjecture based on higher
category theory).  In the present work we state a non-linear version
of this conjecture, and provide an elementary proof of it based on
model category theory and simplicial localization techniques {\em \`a
la} Dwyer-Kan.
\end{blanko*}

\bigskip

The main result of this work is the following. It can reasonably
be considered as a model category version of the well known fact
that the endomorphisms of the unit of a monoidal category form a
commutative monoid.

\begin{theorem*}
  Let $(\M,\tensor,\unit)$ be a monoidal model category. The
simplicial set of derived endomorphisms of the unit,
  $\REnd_{\M}(\unit)$, is a simplicial $2$-monoid (cf.~\ref{mon}).
\end{theorem*}

The theorem applies in particular when $\M$ is the category of
$A$-bimodules over a simplicial monoid $A$.  Then the Hochschild
cohomology $\HH(A)$ is naturally identified with
$\REnd_{\ABim}(A)$,
and hence becomes a $2$-monoid in $\sSet$.  This is what we refer
to as the {\em non-linear analogue of Deligne's conjecture}.

\bigskip

The proof of our theorem relies heavily on ideas of
Segal~\cite{Segal:cat-coh} and
Dwyer-Kan~\cite{Dwyer-Kan:function-complexes}.  Once the
statements have been formulated in terms of Segal categories, the
theorem follows from two easy observations and an application of a
theorem of Dwyer and Kan.

First it is observed that if a monoidal structure on a category is
strictly compatible with a notion of equivalence, the Dwyer-Kan
localization is a {\em monoid\footnote{We warn the reader that the
word \emph{monoid} is used in this work in a much weaker sense
than usual, and always refers to an underlying notion of
\emph{equimorphisms}.  See \ref{mon} for details.} in the category
of simplicial categories} (see \ref{mon}), and by taking the
endomorphism space of the unit object we get instead a monoid in
the category of simplicial monoids, i.e., what we call a
simplicial $2$-monoid.

Second, in the case of a monoidal model category (in the sense of
Hovey~\cite{Hovey:model}), the monoidal operation does not
preserve equivalences, but we observe that Hovey's `unit axiom'
expresses exactly that a suitable equivalent full subcategory has
the strict compatibility and hence we have reduced to the first
case.

The theorem of Dwyer-Kan~\cite{Dwyer-Kan:function-complexes}
describes the derived homomorphisms of a model category (the
simplicial function complexes) in terms of its simplicial
localization.

\bigskip

It is fair to point out that our viewpoint and proof do not seem
to work for the original Deligne conjecture, since currently the
theory of Segal categories does not work well in linear contexts
(like chain complexes), but only in cartesian monoidal contexts.
Also, we did not investigate the relations between $2$-monoids and
simplicial sets with an action of the little $2$-cube operad, and
therefore our version of Deligne's conjecture might be considered
as a bit far from the original one.  However, our original
motivation was not to give an additional proof of Deligne's
conjecture, but rather to try to understand it from a more
conceptual point of view.  The new insight provided by our
approach may also shed light on related subjects.  We also think
it is an interesting application of simplicial localization
techniques.

\bigskip

\begin{blanko*}{Acknowledgements.}
  We are thankful to Andr\'e Hirschowitz and Clemens Berger for
  fruitful discussions, and to V.~Hinich and A.~Voronov for pointing
  out some references on Deligne's conjecture.  The first named author
  wishes to thank the University of Nice for support.
\end{blanko*}

\section{Localization of monoidal coloured categories}

\begin{blanko}{Coloured categories and simplicial localization.}
  By a {\em coloured category} we mean a pair $(C,W)$ where $C$ is a
  category and $W$ is a subclass of arrows, called {\em equimorphisms}
  (or {\em coloured arrows}), closed under composition, and comprising
  all isomorphisms.  Key examples are $\Top$, $\sSet$, and $\Cat$ with
  the usual notions of (weak) equivalences as equimorphisms.  For the
  present purposes, an equally important example is $\sCat$, the
  category of simplicial categories
  (cf.~\cite{Dwyer-Kan:simplicial-localization}): a simplicial functor
  $F:A\to B$ is coloured if $\pi_0 F : \pi_0 A \to \pi_0 B$ is an
  equivalence of categories and for each pair of objects $x,y\in A$,
  the map $A(x,y) \to B(Fx, Fy)$ is a weak equivalence of simplicial
  sets.

  The importance of coloured categories is that they can be localised
  and thus serve as context for expressing weak structures.  Let
  $\CCat$ denote the category whose objects are coloured categories
  and whose arrows are functors that preserve equimorphisms.  The
  classical notion of localization \cite{Gabriel-Zisman} is the
  functor $\Ho: \CCat \to \Cat$ defined by formally inverting all
  equimorphisms.  A much more sophisticated construction is the
  simplicial localization introduced by
  Dwyer-Kan~\cite{Dwyer-Kan:simplicial-localization}, which can be
  seen as a derived version of $\Ho$.  It is a functor $L:\CCat\to
  \sCat$.  It reflects much
  more homotopy theoretic information than the classical localization,
  and in many respects it seems to be the `correct' localization, of
  which the classical localization is just a truncation.  (Indeed, the
  category of connected components of $LC$ is equivalent to $\Ho(C)$.)

  There are several possible ways of turning $\CCat$ into a coloured
  category itself --- the crucial desired property is that $L$ should
  be colour-preserving.  For simplicity we take this as the
  definition: a (colour-preserving) functor $F : (C,W) \to (C',W')$
  between coloured categories is called an {\em equifunctor} if
  $LF$ is an equivalence of simplicial categories.
\end{blanko}

\begin{blanko}{Monoids and $2$-monoids.}\label{mon}
  Let $(S,W)$ be any of the coloured categories mentioned above --- in
  particular $S$ is monoidal with cartesian product as multiplication
  and the singleton object $*$ as unit.  A {\em monoid in $(S,W)$} is a
  functor $X:\Delta\op \to S$ satisfying
  \begin{enumerate}
    \item[{[S0]}] \ $X_0 = *$

    \item[{[S1]}] \ The natural maps $X_k \to X_1 \times \dots \times X_1$
    are equimorphisms ($k\geq 1$).
  \end{enumerate}
  This last axiom is the {\em Segal condition}.  It played a crucial
  r\^ole in Segal's work~\cite{Segal:cat-coh} and was subsequently
  named after him by Tamsamani~\cite{Tamsamani:thesis}.
  A {\em monoid homomorphism} is a natural transformation of such functors.
  The category $\Mon(S)$ of monoids and monoid homomorphisms is
  coloured via the forgetful functor to $S$.  A {\em $2$-monoid in $(S,W)$}
  is by definition a monoid in $\Mon(S)$.

  In the case $S=\sSet$, a monoid is just a Segal category with a
  single object, and a $2$-monoid is a Segal category with a single
  $0$-cell and a single $1$-cell.  For the basic theory of Segal
  categories see \cite{Hirschowitz-Simpson} or \cite{Toen-Vezzosi:0212}.
  Note that this notion of monoid makes sense only
  in cartesian
  monoidal categories (in the ususal sense), since it depends on
  the universal property of the product.
\end{blanko}

\begin{blanko}{Monoidal categories as weak monoids.}
  A monoidal category can be described as a sort of weak monoid
  object in $\Cat$.  The weakness is usually described in terms of $2$-cells
  subject to coherence constraints (e.g., as a bicategory with a
  single object). Here, we will adopt instead the simplicial
  viewpoint, and define
  a monoidal category as a monoid in $\Cat$, in the sense of
  \ref{mon}, conveniently hiding all questions of coherence from the user
  interface.

  This notion is not the same as the usual one defined in terms of
  coherence, but since monoidal categories in either sense are
  equivalent to strict monoidal categories, the two notions lead to
  the same homotopy theory.
  It is not trivial to make specific
  translation between the two languages
  (cf.~Leinster~\cite{Leinster:0002}; see also Segal~\cite{Segal:cat-coh}).
\end{blanko}

\begin{blanko}{Monoidal coloured categories.}\label{mon-col}
  A monoidal coloured category is a monoidal structure on a
  coloured category $(C,W)$ whose structure functors are colour
  preserving.  Precisely,  we define a {\em monoidal coloured
  category} to be a
  functor $\Delta\op\to \CCat$ satisfying [S0] and [S1].
\end{blanko}

\begin{blanko}{Localization of monoidal coloured categories.}\label{Ltensor}
  The way we have set things up it is immediate that the localization
  of a monoidal coloured category is a monoidal simplicial category.
  Indeed, it is just the composite
  $$
  \Delta\op\to \CCat \stackrel{L}{\longrightarrow} \sCat    .
  $$
  Since $L$ preserves the terminal object, preserves products up to
  equivalence (see \cite[Cor.~4.1.2]{Toen-Vezzosi:0212}), and preserves
  equivalences, this composite will again
  satisfy [S0] and [S1].  (Similar remarks hold of course for the
  classical localization.)
\end{blanko}

\begin{blanko}{Endomorphisms of the unit.}
  In fact the simplicial categories appearing in the image are all
  pointed --- the base point is simply the image of $[0]$.  Thus we
  can in a canonical way compose with the endomorphism functor
  $\sCat\lowerstar \to \Mon(\sSet)$, associating to each pointed
  simplicial category the endomorphism monoid of the base point.  This
  is a strict simplicial monoid, and this functor preserves products,
  terminal object, and equimorphisms.  The whole composite is
  therefore a monoid object in $\Mon(\sSet$), i.e., a simplicial
  $2$-monoid.  Hence:
\end{blanko}

\begin{theorem}\label{thm-col}
  Let $((C,W), \tensor, \unit)$ be a monoidal coloured category.
  Then $LC(\unit,\unit)$ is a simplicial $2$-monoid. \qed
\end{theorem}

\section{Localization of monoidal model categories}

Model categories are prominent examples of coloured categories,
and their richer structure allows for important variations on the
localization theme.

\begin{blanko}{The pushout product axiom.}
  Localization of monoidal structure in a model category $M$ was
  considered from the very beginning of model category theory:
  Quillen~\cite{Quillen:homotopical-algebra} observed that in order to
  induce a monoidal structure on $\Ho(\M)$, it is not necessary for a
  monoidal structure on $\M$ to preserve equivalences on the nose, as
  in the general coloured case \ref{mon-col}.
 It is
  enough that the unit is cofibrant and that $\tensor$ satisfies the
  {\em pushout product axiom}: given cofibrations $A_1\to A_2$ and $B_1 \to
  B_2$ then the induced map
  $$
  (A_1\tensor B_2) \coprod_{A_1\tensor B_1} (A_2\tensor  B_1) \longrightarrow A_2
  \tensor B_2
  $$
  is again a cofibration, and if furthermore one of the two original maps
  is a {\em trivial} cofibration then the induced map is too.
  Indeed, in this
  case it follows easily from Ken Brown's Lemma that the full
  subcategory of cofibrant objects $\Mc$ is a monoidal coloured category,
 and in any case $\Mc$ and $\M$ have the same
  homotopy type, so one can induce a monoidal structure on $\Ho(\M)$ by
  taking it from $\Ho(\Mc)$.
\end{blanko}

\begin{blanko}{The unit axiom.}
  Later, Hovey~\cite{Hovey:model} remarked that the requirement that the
  unit be cofibrant can be relaxed: assuming the pushout product axiom
  holds, it is enough that $\M$ satisfies the {\em unit axiom}: for a given
  cofibrant replacement functor $Q: \M \to \Mc$, and for every
  cofibrant $X$, the composite $QI \tensor X \to \unit \tensor X \to X$ is
  an equivalence (and similarly from the right).  In this situation,
  even though the multiplication law on $\Mc$ is not unital, the
  induced multiplication on $\Ho(\M)\simeq \Ho(\Mc)$ does in fact
  acquire a unit.  This justifies the terminology of
  Hovey~\cite{Hovey:model} which has become standard:
\end{blanko}

\begin{definition}
  A {\em monoidal model category} is a model category with a monoidal
  structure satisfying the pushout product axiom and the unit axiom.
\end{definition}

  The following simple observation seems not to have been made before.
  Assume the pushout property axiom holds in $(\M,\tensor,\unit)$, and
  let $\Mci$ denote the full subcategory of all cofibrant objects
  together with the unit.

\begin{lemma}
  The unit axiom holds in $(\M,\tensor,\unit)$ if and only if
  $(\Mci,\tensor,\unit)$ is a monoidal coloured category (i.e., the
  monoidal operation preserves equivalences).
\end{lemma}

\begin{proof}
  Simply note that the unit axiom holds for $Q\unit\isopil \unit$ if
  and only if for {\em any} cofibrant $Z$ with an equivalence $Z \isopil
  \unit$ the conclusion of the unit axiom holds: for cofibrant $X$,
  the map $Z\tensor X \to \unit \tensor X \to X$ is an equivalence.

  Now an equivalence in $\Mci$ is either one between cofibrant objects
  (which case is covered by the pushout product axiom), or of the type
  $Z\isopil \unit$ (the situation just analysed), or $\unit \isopil
  Z$.  But this last type of equivalences is preserved under
  $\tensor$, provided the unit axioms holds, as it readily follows by
  taking a cofibrant replacement of the map and invoking the
  2-out-of-3 axiom for a model category.
\end{proof}

It is easy to see that the monoidal structure induced on $\Ho(\M)$
by Hovey's arguments (resp.~on $L\M$) is merely the one coming
from $\Ho(\Mci)$ (resp.~from $L\Mci$) via the direct construction
of \ref{Ltensor}. One observation is due for this to make sense:

\begin{lemma}\label{l}
  The full embedding $F:\Mci \into \M$ induces an equivalence
  $L\Mci \isopil L\M$ of simplicial categories.
\end{lemma}

\begin{proof}
  In fact this is true for any full subcategory sandwiched between
  $\Mc$ and $\M$. A cofibrant replacement functor $Q:\M\to\Mc\into
  \Mci$ comes with natural transformations $Q\circ F \Rightarrow \id_{\Mci}$
  and $F\circ Q\Rightarrow \id_M$, whose components are equivalences.  By
  standard arguments (see \cite[Lemma 8.1]{Hirschowitz-Simpson} for
  all details), this induces an equivalence after simplicial
  localization.
\end{proof}

\begin{blanko}{Derived endomorphisms.}
  The derived hom set (simplicial function complex) of a pair of
  objects in a model category is usually defined in terms of
  fibrant-cofibrant resolutions functors (see
  e.g.~\cite{Dwyer-Kan:function-complexes}).  We will denote them by
  $\RHom_{\M}(-,-)$.  For two objects $A$ and $B$ in $M$,
  $\RHom_{\M}(A,B)$ is an object in $\sSet$ defined up to
  equivalence.  Of course, $\RHom_{\M}(A,A)$ is denoted by
  $\REnd_{\M}(A)$.

A deep result of Dwyer-Kan~\cite{Dwyer-Kan:function-complexes}
states that this simplicial set  is equivalent to the simplicial
hom sets of the simplicial localization:
$$
\RHom_{\M}(A,B) \simeq L\M(A,B) \simeq L\Mci(A,B) .
$$
(This was actually the original motivation for introducing
simplicial localization.) In particular, by lemma \ref{l} we have
$\REnd_{\M}(\unit) \simeq L\Mci(\unit,\unit)$, and in combination
with Theorem~\ref{thm-col} we get
\end{blanko}

\begin{theorem}\label{thm-model}
  Let $(\M,\tensor,\unit)$ be a monoidal model category. Then
  $\REnd_{\M}(I)$ is a simplicial $2$-monoid.
  \qed
\end{theorem}

Of course, the expression \textit{is} in the above theorem really
means {\em is equivalent to the underlying simplicial set of a
$2$-monoid in $\sSet$}.

\bigskip

\begin{remark}
  In some cases the trick of just adding the non-cofibrant unit by
  hand is not appropriate: for example in K-theory one studies
  Waldhausen categories which are subcategories of the category of
  cofibrant objects, and one cannot just add the unit.  In a similar
  vein, Spitzweck~\cite{Spitzweck} works with a notion of monoidal
  model category with pseudo-unit: this pseudo-unit does not act as a
  unit, but its cofibrant replacements do, up to homotopy.  In these
  cases the important structure is not the `unit' itself but rather
  the space of cofibrant replacements.  These cases are accounted for
  by the theory of monoidal categories with weak units, and more
  generally higher categories with weak identity arrows, where instead
  of strict identities each object has a contractible space of
  up-to-homotopy identity arrows.  The basics of this theory is worked
  out elsewhere; see~\cite{Kock:fairIntro} for an introduction.
  In fact our original approach to the theorem was with weak units,
  but for the present purpose the $\Mci$ trick seems simpler.
\end{remark}

\section{A simplicial version of Deligne's conjecture}

\noindent

\begin{blanko}{Bimodules.}
  Let $A$ be a simplicial monoid (in the strict sense, i.e. a
  simplicial object in the category of monoids), then $A\times A\op$
  is again a simplicial monoid, and we can consider the category of
  $(A\times A\op)$-modules (i.e., simplicial sets with a $(A\times
  A\op)$-action).  $(A\times A\op)$-modules will be called
  $A$-bimodules, and the category of $A$-bimodules is denoted by
  $\ABim$.  This category carries a natural model structure whose
  fibrations and equivalences are induced via the forgetful functor
  $\ABim \to \sSet$ (this is standard, see
  e.g.~Schwede-Shipley~\cite{Schwede-Shipley:9801}).  There is a
  tensor product defined on $\ABim$ as the coequalizer $M\times A \times N
  \rightrightarrows M\times N \to M\tensor_A N$.  The bimodule $A$
  itself is the unit for $\tensor_A$.
\end{blanko}

\begin{lemma}
  $(\ABim, \tensor_A , A)$ is a  monoidal model category.
\end{lemma}

\begin{proof}
  The proof of the lemma relies on a small object argument, using the
  standard generating sets of cofibrations and trivial cofibrations
  (described in \cite{Schwede-Shipley:9801}), as explained in \cite[\S
  4.3]{Hovey:model}.

  Let us recall that the forgetful functor $\ABim \to
  \sSet$ possesses a left adjoint $F : \sSet \to \ABim$,
  sending a simplicial set $X$ to the free $A$-bimodule $F(X)=A\times
  X\times A$.  If $I_{0}$ (resp.~$J_{0}$) is a set of generating
  cofibrations (resp.~trivial cofibrations) in $\sSet$ then
  $I=F(I_{0})$ (resp.~$J=F(J_{0})$) is a set of generating
  cofibrations (resp.~trivial cofibrations) in $\ABim$.

  To prove the pushout product axiom in $\ABim$ it is enough by
  \cite[\S 4.3]{Hovey:model} to notice that for two simplicial sets $X$
  and $Y$ one has a natural isomorphism of $A$-bimodule
  $$
  F(X)\tensor_A F(Y)\simeq F(X\times A\times Y).
  $$
  The pushout product axiom in $\ABim$ is then a direct consequence of
  the well-known facts that the functor $F$ is left Quillen and that the
  pushout product axiom holds in $\sSet$.

  It remains to prove the unit axiom in $\ABim$.  For this we use the
  standard free resolution associated to the forgetful functor $\ABim
  \to \sSet$ (see e.g.~Illusie~\cite{Illusie}).  Let us recall
  that for any $A$-bimodule $M$, one constructs a simplicial object
  $P\lowerstar (M)$ in $\ABim$, together with an augmentation $P_{0}(M)
  \to M$ such that the natural morphism
  $$
  \underset{n\in \Delta\op}{\Hocolim} \; P_{n}(M) \longrightarrow M
  $$
  is an equivalence in $\ABim$.  Furthermore, each $A$-bimodule
  $P_{n}(M)$ is free and given by $P_{n}(M) \df F(P_{n-1}(M))$, and
  the various face and degeneracy morphisms are given by using the
  adjunction between the forgetful functor and $F$.  Since each
  $P_{n}(M)$ is a cofibrant object in $\ABim$ (as it is free), we can
  use Hirschhorn~\cite[Thm.~19.4.2]{Hirschhorn} to see that
  $\Hocolim_{n} P_{n}(M)$ is a cofibrant model for $M$.  To check the
  unit axiom it is therefore enough by \cite[\S 4.3]{Hovey:model} to
  prove that for any simplicial set $X$ the natural morphism
  $$
  (\Hocolim_{n} P_{n}(A))\tensor_A F(X) \longrightarrow
     A\tensor_A F(X)\simeq F(X)
  $$
  is an equivalence. Clearly this morphism is isomorphic to
  $$
  \Hocolim_{n} (P_{n}(A)\tensor_A F(X)) \longrightarrow
    A\tensor_A F(X)\simeq F(X).
  $$
  But, $P_{n}(A)\tensor_A F(X)\simeq P_{n}(A)\times X\times A$, at
  least in $\sSet$, and therefore the morphism is in fact isomorphic,
  as a morphism in $\sSet$, to
  $$
  \Hocolim_{n}(P_{n}(A)\times X\times A) \longrightarrow
    A\times X\times A=F(X).
  $$
  The fact that this last morphism is an equivalence follows simply
  from the fact that $P\lowerstar (A)$ is a simplicial resolution of
  $A$ and that homotopy colimits commute, up to equivalences, with
  products.
\end{proof}

Note that the unit object $A$ of the model category $\ABim$ is not
cofibrant.  Indeed, cofibrant means roughly `free', i.e., direct
sum of copies of $A\times A\op$, but $A$ is rather a quotient.

\begin{blanko}{The Hochschild cohomology.}\label{hh}
  The Hochschild cohomology of a simplicial monoid $A$ can naturally
  be defined as
  $$
  \HH(A) \df \REnd_{\ABim}(A).
  $$
  The Hochschild cohomology of a simplicial monoid is clearly a
  homotopy version of its centre.  Indeed, if $M$ is a monoid (in the
  category of sets), the endomorphisms of $M$ as a $M$-bimodule is
  naturally isomorphic to the centre of $M$.

  There exist also more explicit descriptions of the Hochschild
  cohomology of a simplicial monoid $A$, more in the style of the
  Hochschild complex of an associative algebra.  They can be obtained
  by taking explicit cofibrant replacement of $A$ as an $A$-bimodule.
\end{blanko}

\bigskip

Finally, Theorem~\ref{thm-model} applies, yielding the following
corollary, which we call a theorem for emphasis:

\begin{theorem}\label{deligne-conj}
  Let $A$ be a simplicial monoid.  Then the Hochschild cohomology
  $\HH(A)$ is a simplicial $2$-monoid.
  \qed
\end{theorem}

\section{Higher dimensional generalization}

Theorem \ref{deligne-conj} can be generalized in the following way
in terms of Segal categories (starting from the observation that a
simplicial monoid is a Segal $1$-monoid).  First of all, the
definition of monoids in a coloured category as described in
\ref{mon} can be iterated.  Starting by letting $\kat{0-SeMon}$ be
the coloured category of simplicial sets and equivalences, one
defines (for $d\geq 1$) the coloured category $\kat{d-SeMon}$ of
{\em Segal $d$-monoids} as
$$
\kat{d-SeMon} \df
\Mon\big(\textbf{(}\kat{d-1}\textbf{)}\kat{-SeMon}\big),
$$
in the sense of \ref{mon}.

Any Segal $d$-monoid $M$ has an underlying Segal $1$-monoid, which
up to equivalence can be chosen to be a simplicial monoid in the
usual sense (i.e.~a simplicial object in the category of monoids),
cf.~e.g.~\cite[\S 8]{Hirschowitz-Simpson}.  We define {\em the
Hochschild cohomology of a Segal $d$-monoid} to be the Hochschild
cohomology of its underlying simplicial monoid, as defined in
\ref{hh}.  Theorem \ref{deligne-conj} now has the following
generalization.

\begin{theorem}\label{d-deligne-conj}
  Let $A$ be a Segal $d$-monoid.  Then the Hochschild cohomology
  $\HH(A)$ is a Segal $(d+1)$-monoid.
  \qed
\end{theorem}

We will not include the proof of this theorem as it uses the
theory of Segal categories and the so-called
\textit{strictification theorem} stated in
\cite{Toen-Vezzosi:0212}.  It would be interesting however to have
a model category proof of Theorem~\ref{d-deligne-conj}.  A
possible approach would be through a suitable notion of
\textit{iterated model category}, which roughly would be an
iterated monoidal category in the sense of
\cite{Balteanu-et.al.:9808}, together with a compatible model
category structure.  Our Theorem \ref{thm-model} should then
generalize as follows: if $M$ is a $d$-times iterated monoidal
model category then $\REnd_\M(\unit)$ is a $(d+1)$-monoid.


\begin{thebibliography}{10}

\bibitem{Balteanu-et.al.:9808}
{\sc Cornel B\u alteanu, Zbigniew Fiedorowicz, Roland
Schw{\"a}nzl, {\rm and }Rainer Vogt}.
\newblock {\em Iterated monoidal categories}.
\newblock   Preprint, math.AT/9808082.

\bibitem{Batanin:0207}
{\sc Michael Batanin}.
\newblock {\em The Eckmann-Hilton argument,
higher operads and $E_{n}$-spaces}.
\newblock   Preprint, math.CT/0207281.

\bibitem{Berger-Fresse:0109}
{\sc Clemens Berger {\rm and }Benoit Fresse}.
\newblock {\em Combinatorial operad actions on cochains}.
\newblock Preprint, math.AT/0109158.

\bibitem{Dwyer-Kan:simplicial-localization}
{\sc William~G. Dwyer {\rm and }Daniel~M. Kan}.
\newblock {\em Simplicial localizations of categories}.
\newblock J.~Pure Appl. Algebra {\bf 17} (1980), 267--284.

\bibitem{Dwyer-Kan:function-complexes}
{\sc William~G. Dwyer {\rm and }Daniel~M. Kan}.
\newblock {\em Function complexes in homotopical algebra}.
\newblock Topology {\bf 19} (1980), 427--440.

\bibitem{Gabriel-Zisman}
{\sc Peter Gabriel {\rm and }Michel Zisman}.
\newblock {\em Calculus of fractions and homotopy theory}.
\newblock Ergebnisse der Mathematik und ihrer Grenzgebiete, Band 35.
  Springer-Verlag New York, Inc., New York, 1967.

\bibitem{Gerstenhaber-Voronov:9409}
{\sc Murray Gerstenhaber {\rm and }Alexander A. Voronov}.
\newblock {\em Homotopy {$G$}-algebras and moduli space operad},
\newblock Internat. Math. Res. Notices {\bf 3} (1995) 141--153 (electronic).

\bibitem{Getzler-Jones:9403}
{\sc Ezra Getzler {\rm and }John D.~S. Jones}.
\newblock {\em Operads, homotopy algebra and iterated integrals for double
  loop spaces}.
\newblock Preprint, hep-th/9403055.

\bibitem{Hirschhorn}
{\sc Philip S. Hirschhorn}.
\newblock {\em Model categories and their localizations},
  vol.~99 of Mathematical Surveys and Monographs.
\newblock American Mathematical Society, Providence, RI, 2003.
\newblock (A preliminary version was widely circulated under the
  title {\em Localization of model categories}. This is the version
  we actually refer to.)

\bibitem{Hirschowitz-Simpson}
{\sc Andr{\'e} Hirschowitz {\rm and }Carlos Simpson}.
\newblock {\em Descente pour les {$n$}-champs}.
\newblock Preprint, math.AG/9807049.

\bibitem{Hovey:model}
{\sc Mark Hovey}.
\newblock {\em Model categories}, vol.~63 of Mathematical Surveys and
  Monographs.
\newblock American Mathematical Society, Providence, RI, 1999.

\bibitem{Hu-Kriz-Voronov}
{\sc Po Hu, Igor Kriz, {\rm and }Alexander A. Voronov}.
{\em On Kontsevich's Hochschild cohomology conjecture}.
\newblock Preprint, available at http://www.math.umn.edu/{\~{ }}voronov/igor.ps.


\bibitem{Illusie}
{\sc Luc Illusie}.
\newblock {\em Complexe cotangent et d\'eformations. {I}}.
\newblock No.~239 in Lecture Notes in Mathematics.
Springer-Verlag, Berlin, 1971.

\bibitem{Kock:fairIntro}
{\sc Joachim Kock}.
\newblock {\em A notion of weak identity arrows in higher categories}.
\newblock Manuscript, available at
  http://www-math.unice.fr/{\~{}}kock/cats/aarhus.pdf.

\bibitem{Kontsevich:9904}
{\sc Maxim Kontsevich}.
\newblock {\em Operads and motives in deformation quantization}.
\newblock Lett. Math. Phys.  {\bf 48} (1999), 35--72.
\newblock (math.QA/9904055).

\bibitem{Kontsevich-Soibelman:0001}
{\sc Maxim Kontsevich {\rm and }Yan Soibelman}.
\newblock {\em Deformations of algebras over operads and the {D}eligne
  conjecture}.
\newblock In {\em Conf\'erence Mosh\'e Flato 1999, Vol. I (Dijon)}, vol.~21 of
  Math. Phys. Stud., pp. 255--307. Kluwer Acad. Publ., Dordrecht, 2000.
\newblock (math.QA/0001151).

\bibitem{Leinster:0002}
{\sc Tom Leinster}.
\newblock {\em Homotopy algebras for operads}.
\newblock Preprint, math.CT/0002180.

\bibitem{McClure-Smith:9910}
{\sc James~E. McClure {\rm and }Jeffrey~H. Smith}.
\newblock {\em A solution of {D}eligne's {H}ochschild cohomology conjecture}.
\newblock In {\em Recent progress in homotopy theory (Baltimore, MD, 2000)},
  vol. 293 of Contemp. Math., pp. 153--193. Amer. Math. Soc., Providence, RI,
  2002.
\newblock (math.QA/9910126).

\bibitem{Quillen:homotopical-algebra}
{\sc Dan Quillen}.
\newblock {\em Homotopical algebra}.
\newblock No.~43 in Lecture Notes in Mathematics. Springer-Verlag, Berlin,
  1967.

\bibitem{Schwede-Shipley:9801}
{\sc Stefan Schwede {\rm and }Brooke Shipley}.
\newblock {\em Algebras and modules in monoidal model categories}.
\newblock Proc. London Math. Soc. (3) {\bf 80} (2000), 491--511.
\newblock (math.AT/9801082.)

\bibitem{Segal:cat-coh}
{\sc Graeme Segal}.
\newblock {\em Categories and cohomology theories}.
\newblock Topology {\bf 13} (1974), 293--312.

\bibitem{Spitzweck}
{\sc Markus Spitzweck}.
\newblock {\em Operads, algebras and modules in model categories and motives}.
\newblock Preprint, 2001.

\bibitem{Tamarkin:9803}
{\sc Dmitry E. Tamarkin}.
\newblock {\em Another proof of {M}.~{K}ontsevich formality 
   theorem for {$\mathbb{R}^n$}}.
\newblock Preprint, math.QA/9803025.

\bibitem{Tamsamani:thesis}
{\sc Zouhair Tamsamani}.
\newblock {\em Sur des notions de $n$-cat\'egorie et $n$-groupo\"\i de non
  strictes via des ensembles multi-simpliciaux}.
\newblock $K$-Theory {\bf 16} (1999), 51--99.
\newblock (alg-geom/9512006 and alg-geom/9607010.)

\bibitem{Toen-Vezzosi:0212}
{\sc Bertrand To{\"e}n {\rm and }Gabriele Vezzosi}.
\newblock {\em Segal topoi and stacks over {S}egal categories}.
\newblock Preprint, math.AG/0212330.

\bibitem{Voronov:9908}
{\sc Alexander A. Voronov}.
\newblock {\em Homotopy Gerstenhaber algebras}. 
\newblock In {\em Conf\'erence Mosh\'e Flato 1999, Vol. II (Dijon)}, vol.~22 of
Math. Phys. Stud., pp. 307--331,  Kluwer Acad.
Publ., Dordrecht, 2000. 
\newblock (math.QA/9908040.)

\end{thebibliography}

\end{document}